\documentclass[11pt]{article}
\usepackage[english]{style}
\usepackage{amsmath,amsfonts,graphics,xypic}

\xyoption{all}

\DeclareFontFamily{U}{rsf}{}
\DeclareFontShape{U}{rsf}{m}{n}{
  <5> <6> rsfs5 <7> <8> <9> rsfs7 <10-> rsfs10}{}
\DeclareMathAlphabet{\mathscr}{U}{rsf}{m}{n}

\DeclareMathAlphabet{\mathgth}{U}{euf}{m}{n}

\DeclareFontFamily{U}{cyr}{}
\DeclareFontShape{U}{cyr}{m}{n}{
  <5> wncyr5 <6> wncyr6 <7> wncyr7 <8> wncyr8 <9> wncyr9 <10-> wncyr10}{}
\DeclareMathAlphabet{\mathcyr}{U}{cyr}{m}{n}

\input cyracc.def

\makeatletter
\def\operator@font{\sf}
\makeatother

\newcommand{\gh}{\mathgth{h}}

\newcommand{\ccH}{{\mathscr H}}

\newcommand{\cO}{{\mathscr O}}

\newcommand{\cZ}{{\mathscr Z}}

\newcommand{\bi}{{\bar{\imath}}}
\newcommand{\bj}{{\bar{\jmath}}}

\newcommand{\HKR}{{\mathsf{HKR}}}

\newcommand{\D}{{\mathbf D}}

\newcommand{\chk}{{\scriptscriptstyle\vee}}
\newcommand{\R}{\mathbf{R}}

\newcommand{\Ld}{\mathbf{L}}

\newcommand{\bbk}{\mathbf{k}}

\newcommand{\bbE}{\mathbb{E}}
\newcommand{\bbF}{\mathbb{F}}
\newcommand{\bbN}{\mathbb{N}}
\newcommand{\bbT}{\mathbb{T}}

\renewcommand{\S}{{\mathsf{Sym}}}

\newcommand{\Hom}{{\mathsf{Hom}}}

\DeclareMathOperator{\Tor}{Tor}

\DeclareMathOperator{\Ext}{Ext}

\newcommand{\T}{\mathbb{T}}

\newcommand{\ra}{\rightarrow}

\newcommand{\lra}{\longrightarrow}

\newcommand{\Z}{\mathbf{Z}}

\newcommand{\DSch}{\mathgth{DSch}}

\newcommand{\iso}{\cong}

\newcommand{\HH}{H\mskip-2mu H}

\renewcommand{\phi}{\varphi}

\author{%
Dima Arinkin\thanks{Mathematics Department,
University of Wisconsin--Madison, 480 Lincoln Drive, Madison, WI
53706, USA, {\em e-mail: }{\tt arinkin@math.wisc.edu,
andreic@math.wisc.edu}},\ \ 
Andrei C\u ald\u araru,\footnotemark[1]\ \  M\'arton
Hablicsek\thanks{Mathematics Department, University of Pennsylvania,
David Rittenhouse Lab, 209 S.\ 33rd Street, Philadelphia, PA 19104,
USA, {\em e-mail: }{\tt mhabli@math.upenn.edu}}
}

\title{Formality of derived intersections and the orbifold HKR isomorphism}

\date{}

\begin{document}

\maketitle

\begin{abstract}
  We study when the derived intersection of two smooth subvarieties of
  a smooth variety is formal.  As a consequence we obtain a derived
  base change theorem for non-transversal intersections.  We also
  obtain applications to the study of the derived fixed locus of a
  finite group action and argue that for a global quotient orbifold
  the exponential map is an isomorphism between the Lie algebra of the
  free loop space and the loop space itself.  This allows us to give
  new proofs of the HKR decomposition of orbifold Hochschild
  (co)homology into twisted sectors.
\end{abstract}

\section{Introduction}

\paragraph
\label{subsec:diagintro}
Let $S$ be a smooth variety over a field of characteristic zero
and let $X$ and $Y$ be smooth subvarieties of $S$.  We shall assume
that $X$ and $Y$ intersect cleanly (meaning that their scheme
theoretic intersection $W = X \times_S Y$ is smooth) but not
necessarily transversely.  Derived algebraic geometry associates to
this data a geometric object, the {\em derived intersection} of $X$
and $Y$,
\[ W' = X\times^R_S Y. \]
It is a differential graded (dg) scheme whose structure complex is
constructed by taking the derived tensor product of the structure
sheaves of $X$ and $Y$.  (The reader unfamiliar with the subject of dg
schemes is referred to Section~\ref{sec:dgsch}.) The underived
intersection $W$ naturally sits inside $W'$ as a closed subscheme.

\noindent
We organize these spaces and the maps between
them in the diagram
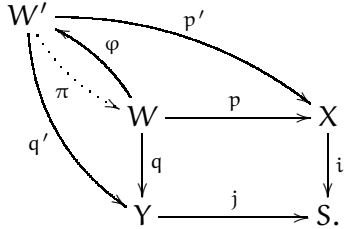
\[ 
\xymatrix{W' \ar@/_/@{.>}[dr]_{\pi} \ar@/_1pc/[ddr]_{q'}
  \ar@/^1pc/[drrr]^{p'}& & & \\ & W\ar@/_/[ul]_-{\phi}
  \ar[rr]^p\ar[d]^-q & & X\ar[d]^-i\\ & Y\ar[rr]^j & & S.\\}
\]
The purpose of this note is to understand when $W'$ is as simple as
possible.  Our main result (Theorem~\ref{thm:mainthm}) makes this
precise in two ways.  In algebraic terms it describes when $W'$ is
{\em formal} in the sense of~\cite{DelGriMorSul}.  In geometric terms
it gives a necessary and sufficient condition for the existence of a
map $\pi:W'\ra W$ exhibiting $W'$ as the total space of a (shift of a)
vector bundle over $W$.  When this holds we gain a geometric
understanding of the structure of the maps $\phi$, $p'$ and $q'$:
$\phi$ is the inclusion of the zero section of the bundle, and $p'$
and $q'$ factor through the bundle map $\pi$.

\paragraph
The problem we study originates in classical intersection theory.  While the
scheme-theoretic intersection $W$ is determined algebraically by the
{\em underived} tensor product
\[ \cO_W = \cO_X \otimes_{\cO_S} \cO_Y, \] 
Serre~\cite{Ser} argued that in order to obtain a theory with good
formal properties we need to use instead the {\em derived} tensor
product
\[ \cO_{W'} = \cO_X \otimes^L_{\cO_S} \cO_Y. \] 
Since $\cO_{W'}$ is naturally a commutative dg
algebra we can regard it as the structure complex of a dg scheme $W'$.

\paragraph
\label{subsec:calcTor}
For classical applications it suffices to work with the class of
$\cO_{W'}$ in K-theory.  Put differently, we only need to know the
sheaves
\[ \ccH^{-*}(\cO_{W'}) = \Tor_*^{S}(\cO_X, \cO_Y). \] 
A local calculation as in~\cite[Proposition A.3]{CalKatSha} shows that
these sheaves are the exterior powers $\wedge^* E^\chk$ of the {\em
  excess bundle} $E$, the vector bundle on $W$ defined by
\[ E = \frac{T_S}{T_X+T_Y}. \]
(We omit writing the restrictions from $X$, $Y$, $S$ to $W$ in
formulas like the one above.  So when we write $T_X$ we mean $T_X|_W$,
the restriction of the tangent bundle of $X$ to $W$.) The excess
bundle $E$ is a vector bundle on $W$ of rank $\dim S + \dim W - \dim X
-\dim Y$ which measures the failure of the intersection to be
transversal.

\paragraph
\label{subsec:extss}
For certain problems it is not enough to know just the
$\Tor$ sheaves $\ccH^k(\cO_{W'})$; we need to understand the full dg
algebra $\cO_{W'}$.  For example there is considerable interest in
computing $\Ext^*_S(i_*F, j_* G)$ for vector bundles $F$ and $G$ on
$X$ and on $Y$.  These groups can be computed using the spectral
sequence
\[ {}^2E^{pq} = H^p(W, F^\chk \otimes G \otimes \omega_{W/Y}\otimes \wedge^{q-m} E)
\Rightarrow \Ext^{p+q}_S(i_* F, j_* G), \] 
from~\cite[Theorem A.1]{CalKatSha}, where $\omega_{W/Y}$ denotes the relative dualizing sheaf of the embedding $W\ra Y$ and $m$ denotes the codimension of $W$ inside $Y$.  (Again, we omit the restrictions
of $F$ and $G$ to $W$.)  The differentials in this spectral sequence
arise as obstructions to splitting the canonical filtration on
$\cO_{W'}$, that is, they vanish if there is an isomorphism
\[ \cO_{W'} \iso \bigoplus_k \ccH^k(\cO_{W'})[-k]. \]

\paragraph
In the above discussion we have skated over an important detail.  The
splitting of $\cO_{W'}$ is not an intrinsic property of the dg scheme
$W'$; rather, the concept only makes sense for a morphism from $W'$ to
to a base scheme.  We discuss its relation to a more general notion,
that of formality.

Consider a dg scheme $Z'$ which is affine over an ordinary scheme $Z$,
i.e., the dg scheme $Z'$ is endowed with a structure morphism $s:Z'
\ra Z$ such that $s$ is affine.  We shall consider two related dg
algebras over $\cO_Z$.  One is $s_*\cO_{Z'}$; the other is its
associated graded counterpart
\[ \cO_{\hat{Z}'} = \bigoplus_k \ccH^k(s_*\cO_{Z'})[-k]. \] 
Note that the right hand side inherits an associative commutative
product structure from that of $\cO_{Z'}$, so it can be regarded as
the structure complex of a dg scheme $\hat{Z}'$ over $Z$.

We shall say that $Z'$ is formal over $Z$ if there exists an
isomorphism $Z' \iso \hat{Z}'$ in the category of derived dg schemes
over $Z$.  (See Section~\ref{sec:dgsch} for further details.)  This is
equivalent to saying that the dg algebra $s_* \cO_{Z'}$ is a formal
$\cO_Z$-dg algebra, that is, there exists an isomorphism of {\em dg
  algebras} over $\cO_Z$
\[ s_*\cO_{Z'} \iso \cO_{\hat{Z}'}. \]
(The terminology is inspired by the work~\cite{DelGriMorSul}
of Deligne, Griffiths, Morgan, and Sullivan, where the authors define
a smooth manifold to be formal if its de Rham dg algebra is
formal in the above sense.)

Note that in particular $Z'$ being formal over $Z$ implies that the
complex $s_*\cO_{Z'}$ is split in $\D(Z)$ (it is isomorphic, as an
$\cO_Z$-module, to the sum of its cohomology sheaves).  In our context
this is a non-trivial condition, unlike the situation
in~\cite{DelGriMorSul}.

\paragraph
The derived intersection $W' = X\times^R_S Y$ can be viewed as a dg
scheme over several base schemes: either one of $X$, $Y$, $S$, or
$X\times Y$ can naturally serve as an underlying scheme for $W'$.
(However, note that in general we can not present $W'$ as a dg scheme
over $W$.)  Our primary motivation for studying derived intersections
comes from our desire to understand the degeneration of the spectral
sequence in~(\ref{subsec:extss}).  For this purpose it is most useful
to regard $W'$ as a dg scheme over $X\times Y$.  Indeed, in this
approach the structure sheaf $\cO_{W'}$ of $W'$ is the kernel of the
functor $j^* i_*:\D(X) \ra \D(Y)$ between dg enhancements $\D(X)$,
$\D(Y)$ of the derived categories of $X$ and $Y$.  (We omit the $\R$'s
and $\Ld$'s in front of derived functors for simplicity.)  

\paragraph
Formality of $W'/X\times Y$ turns out to be closely related to
properties of the inclusion $W\hookrightarrow W'$.  We shall say that
a map of spaces $W \ra W'$ {\em splits} if it admits a left inverse.
(The term ``splitting'' may be misleading: it might help the reader to
think of split embeddings as algebro-geometric analogues of
deformation retractions in topology.)  If $W \ra W'$ is a closed
embedding we shall say that it {\em splits to first order} if the
induced map $W \ra W'{}^{(1)}$ splits, where $W'{}^{(1)}$ is the first
infinitesimal neighborhood of $W$ inside $W'$.

The above concepts also make perfect sense for spaces (schemes, dg
schemes) over a fixed base scheme, in which case we require the
inverse map to be a map over the base scheme.  
\medskip

\noindent
We are now ready to state the main theorem of the paper
 \footnote{While working on the final draft of this paper
    the authors became aware that a closely related result was obtained
    independently and at about the same time by Grivaux~\cite{Gri}.}.

\begin{Theorem}
\label{thm:mainthm}
The following statements are equivalent.
\begin{itemize}
\item[(1)] There exists an isomorphism of dg functors $\D(X) \ra \D(Y)$
\[ j^* i_* (\,-\,) \iso q_*(p^*(\, - \,) \otimes \S(E^\chk[1])). \]
\item[(2)] There exists an isomorphism $W'\iso \bbE[-1]$ of dg
  schemes over $X\times Y$.
\item[(3)] $W'$ is formal as a dg scheme over $X\times Y$.
\item[(4)] The inclusion $W \ra W'$ splits over $X\times Y$.
\item[(5)] The inclusion $W \ra W'$ splits to first order over $X\times Y$.
\item[(6)] The short exact sequence
\[ 0 \ra T_X + T_Y \ra T_S \ra E \ra 0 \]
of vector bundles on $W$ splits.
\end{itemize}
\end{Theorem}

\paragraph
The above theorem can be seen as a generalization of several classical
results: base change for flat morphisms or, more generally, for
$\Tor$-independent morphisms; the Hochschild-Kostant-Rosenberg
isomorphism for schemes~\cite{Swa}; and the formality theorem for
derived self-intersections of the first two authors~\cite{AriCal}.  A
slightly modified version of this theorem was used in a twisted
context in~\cite{AriCalHab} to prove a theorem on the formality of the
twisted de Rham complex.

\paragraph
\label{subsec:derfixloc}
However, the main application we have in mind for
Theorem~\ref{thm:mainthm} is in the study of derived fixed loci.  Let
$\phi$ be a finite-order automorphism of a smooth variety $Z$.  We are
interested in the fixed locus $W$ of $\phi$,
\[ W= Z^\phi = \{ z \in Z~|~ \phi(z) = z \}. \]
This fixed locus can be studied using intersection theory, as we can
view $W$ as the intersection (inside $Z\times Z$) of the diagonal
$\Delta$ and the graph $\Delta^\phi$ of $\phi$,
\[ W = \Delta \times_{Z\times Z} \Delta^\phi. \]

\paragraph
This description makes it clear that the expected dimension of the
fixed locus is zero.  Whenever $W$ is positive dimensional the cause
is a failure of transversality of $\Delta$ and $\Delta^\phi$.  It then
makes sense to study the {\em derived fixed locus} of $\phi$, $W'$,
which we define as the derived intersection
\[ W' = \Delta\times^R_{Z\times Z} \Delta^\phi. \] 
The excess intersection bundle $E$ for this problem is easily seen to
be precisely $T_W$, the tangent bundle of the underived fixed locus
$W$. 
\medskip

\noindent
In this setup Theorem~\ref{thm:mainthm} allows us to get the
following geometric characterization of the derived fixed locus $W'$.

\begin{Corollary}
\label{cor:fpl}
The derived fixed locus $W'$ is isomorphic, as a dg scheme over
$Z\times Z$, to the total space over $W$ of the dg vector bundle
$T_W[-1]$,
\[ W' \iso \bbT_W[-1]. \]
\end{Corollary}

\paragraph
We apply the above result to the study of loop spaces of orbifolds.
Recall that for a space $X$ one defines the {\em free loop space} of
$X$ as 
\[ LX = X \times^R_{X\times X} X. \]
It is naturally a formal derived group scheme over $X$.  

When $X$ is a smooth scheme, the relative Lie algebra of $LX/X$ can be
identified with $\T_X[-1]$, the total space of the shifted tangent
bundle $T_X[-1]$.  In this case the exponential map is an isomorphism
(commonly called the Hochschild-Kostant-Rosenberg isomorphism, or HKR)
\[ \exp: \T_X[-1] \stackrel{\sim}{\lra} LX. \] 
The non-trivial part of the above statement is the fact that the
exponential map is an isomorphism not only in a formal neighborhood of
the origin, but in fact extends to the whole group.  (This follows
from the fact that the loop space $LX$ is a nilpotent extension of $X$.)

\paragraph 
The above statement is known to fail for Artin stacks.  See for
example Ben-Zvi and Nadler~\cite{BenNad}, where the
authors prove that in this case only the formal version of the HKR
isomorphism holds.

By contrast, as an application of the formality of derived fixed loci,
we prove that the HKR isomorphism theorem still holds for global
quotient orbifolds (global quotient Deligne-Mumford stacks).  
The setting we will work with is as follows. Let $G$ be a finite
group acting on a smooth variety $Z$, and let $\cZ$ be the quotient
stack $[Z/G]$.  Denote by $I\cZ$ the (underived) inertia stack
\[ I\cZ = \cZ \times_{\cZ \times \cZ} \cZ, \] 
by $\T_{I\cZ}[-1]$ its shifted tangent bundle, and by $L\cZ$ 
the free loop space of $\cZ$
\[ L\cZ = \cZ \times_{\cZ \times \cZ}^R \cZ. \]

\paragraph
{\bf Theorem (Orbifold HKR isomorphism).}
\label{thm:orb1}
{\em Let $\cZ$ be a global quotient orbifold.  Then there exists a
  canonical isomorphism 
  \[ \exp: \T_{I\cZ}[-1] \stackrel{\sim}{\lra} L\cZ \] 
between the shifted tangent bundle $\T_{I\cZ}[-1]$ of its inertia orbifold and its
free loop space $L\cZ$.}

\paragraph
As in the case of the usual HKR isomorphism for smooth schemes, the
above theorem allows us to give a decomposition of the Hochschild
(co)homology of the orbifold $\cZ$.  In order to state this
decomposition in more concrete terms we need the following notations
for $g\in G$:
\begin{itemize}
\item[--] $Z^g$ is the (underived) fixed locus of $g$ in $Z$;
\item[--] $i_g$ is the closed embedding of $Z^g$ in $Z$;
\item[--] $c_g$ is the codimension of $Z^g$ in $Z$;
\item[--]$\omega_g$ is the relative dualizing bundle of the embedding
  $i_g$, that is, the top exterior power of the normal bundle
  $N_{Z^g/Z}$ of $Z^g$ in $Z$,
  \[ \omega_g = \wedge^{c_g} N_{Z^g/Z}; \]
\item[--] $T_g$ is the vector bundle on $Z^g$ obtained by taking
  coinvariants of $T_Z|_{Z^g}$ with respect to the action of $g$;
\item[--] $\Omega_g^j$ is the dual, along $Z^g$, of $\wedge^j T_g$;
\item[--] $\S(\Omega_g^1[1])$ is the symmetric algebra of
  $\Omega_g^1[1])$, i.e., the object of $\D(Z^g)$
\[ \S(\Omega_g^1[1]) = \oplus \Omega_g^j[j]. \]
\end{itemize}
\vspace*{-2mm}

\noindent
With these notations we can phrase the following consequence of
Theorem~\ref{thm:orb1}. 

\begin{Corollary}
\label{cor:orb2}
The two projections $p', q':L\cZ \ra \cZ'$ are homotopic (and hence
equal in the derived category of dg stacks).  There are natural
isomorphisms of dg functors $\D(\cZ) \ra \D(\cZ)$
\[ \Delta^* \Delta_* (\,-\,) \iso q'_* p'{}^* (\,-\,) \iso \,-\,
\otimes q'_* \cO_{L\cZ}. \]  
The object $q'_* \cO_{L\cZ} \in \D(\cZ)$ is represented by the
$G$-equivariant object of $\D(Z)$
\[ \bigoplus_{g\in G} i_{g,*} \,\S(\Omega_g^1[1]). \]
\newpage
\noindent
Therefore 
\begin{itemize}
\item[(1)] $\Delta^*\Delta_* \cO_\cZ = \bigoplus_{g\in G} i_{g,*}
  \S(\Omega_Z^g[1])$.
  \item[(2)] $\HH_*(\cZ) = \left ( \bigoplus_{g\in G} \bigoplus_{q-p =
        *} H^p(Z^g, \Omega_g^q)  \right )^G.$
  \item[(3)] $HH^*(\cZ) = \left ( \bigoplus_{g \in G} \bigoplus_{p+q =
        *} H^{p-c_g}(Z^g, \wedge^q T_g \otimes \omega_g )\right )_G$.
\end{itemize}
\end{Corollary}
\medskip

\noindent
This result generalizes known Hochschild-Kostant-Rosenberg isomorphisms
for orbifolds, for example those of Baranovsky~\cite{Bar} and
Ganter~\cite{Gan}.

\paragraph
The paper is organized as follows.  In Section~\ref{sec:dgsch} we
collect some general results about dg schemes in the sense of
Ciocan-Fontanine and Kapranov.  In particular we discuss the concept
of dg schemes relative to a base scheme and the concept of formality.
We construct presentations of the derived intersection $W'$ over $X$,
$Y$, $X\times Y$, and $S$.  In Section~\ref{sec:mainthm} we present
the proof of Theorem~\ref{thm:mainthm}.  In the final section of the
paper we discuss applications to orbifolds, and present proofs of
Corollary~\ref{cor:fpl}, Theorem~\ref{thm:orb1}, and
Corollary~\ref{cor:orb2}.

\paragraph
\textbf{Conventions.}  
We work over a field of characteristic zero.  The same results also
hold when the characteristic of the ground field is sufficiently
large; we shall make it explicit in the statement of each theorem how
large the characteristic needs to be for the results to hold.  All
schemes are assumed to be smooth, quasi-projective over this field.

\paragraph
\textbf{Acknowledgments.}  The present project originates in an old
discussion the second author had around 1996 with Dan Abramovich.  We
have benefited from stimulating conversations with Tony Pantev.  The
authors are supported by the National Science Foundation under Grants
No.\ DMS-0901224, DMS-1101558, and DMS-1200721.

\section{Background on dg schemes}
\label{sec:dgsch}

In this section we review some facts from the basic theory of
differential graded schemes, following the work of Ciocan-Fontanine
and Kapranov~\cite{CioKap}.  We emphasize the point of view that a dg
scheme $Z' = (Z, \cO_{Z'})$ should be thought of as a dg scheme {\em
  over} $Z$, and explain how the derived intersection $W' =
X\times_S^R Y$ can be viewed in a natural way as a dg scheme over $X$,
$Y$, $X\times Y$, or $S$.

\paragraph
Following Ciocan-Fontanine and Kapranov~\cite{CioKap}, a {\em
  differential graded scheme} $Z'$ is a pair $(Z, \cO_{Z'})$
consisting of an ordinary scheme $Z$, the {\em base scheme} of $Z'$,
and a complex of quasi-coherent sheaves $\cO_{Z'}^\cdot$ on $Z$, the
{\em structure complex} of $\Z'$.  The complex $\cO_{Z'}$ is assumed
to be endowed with the structure of a commutative dg algebra over
$\cO_Z$, and must satisfy
\begin{itemize}
\item[1.] $\cO_{Z'}^i = 0$ for $i>0$;
\item[2.] $\cO_{Z'}^0 = \cO_Z$.
\end{itemize}
Maps between dg schemes are obtained by a localization procedure
similar to the one that leads to the construction of derived
categories. In a first stage morphisms of dg schemes are considered as
maps of ringed spaces.  For dg schemes $Z' = (Z, \cO_{Z'})$ and $W' =
(W, \cO_{W'})$ a morphism $Z' \ra W'$ consists of a map of schemes
$f:Z\ra W$ along with a map of dg algebras $f^\#:f^*\cO_{W'} \ra
\cO_{Z'}$.  In the resulting category we have a natural notion of
quasi-isomorphisms of dg schemes -- those morphisms $(f, f^\#)$ for
which $f^\#$ is a quasi-isomorphism of complexes of sheaves.  Formally
inverting those quasi-isomorphisms produces a category $\DSch$, the
right derived category of schemes.

\paragraph
Because quasi-isomorphisms become isomorphisms in $\DSch$, isomorphic
dg schemes can be presented over different base schemes.  Thus the
base scheme is not an intrinsic part of a dg scheme in $\DSch$.  For
certain purposes, however, it is useful to be able to refer to the
base scheme of a dg scheme.  Instead of carrying over this additional
data, we give an alternative way of looking at the relationship
between a dg scheme $Z'$ and its supporting scheme $Z$.

The definition of dg schemes implies that the structure complex
$\cO_{Z'}$ of a dg scheme $Z' = (Z, \cO_{Z'})$ admits a natural
morphism of dg algebras $\cO_Z \ra \cO_{Z'}$ (where $\cO_Z$ is
regarded as a complex concentrated in degree zero). This shows that a
dg scheme $Z'$ presented over a base scheme $Z$ comes with a canonical
morphism $Z'\ra Z$.

\paragraph
This observation motivates us to study dg schemes {\em over a fixed
  scheme} $Z$ instead of arbitrary dg schemes.  These are dg schemes
$Z'$ endowed with a morphism $Z' \ra Z$.  (We shall mostly be
concerned with the situation when this morphism is {\em affine} --
this is the case when the dg scheme $Z'$ is presented over $Z$.  But
the concept makes sense in general.)  As in the theory of schemes,
morphisms of dg schemes over $Z$ are morphisms between dg schemes
which commute with the structure morphisms.

\paragraph
We now turn to discussing the construction of derived intersections
over various bases.  We place ourselves in the context described in
the introduction, with $X$ and $Y$ subschemes of $S$.  The structure
complex of the derived intersection $W' = X\times_S^R Y$ is obtained
by taking the derived tensor product $\cO_{W'} = \cO_X
\otimes_{\cO_S}^\Ld\cO_Y$.

The main question we want to address is over what base scheme should
the complex $\cO_{W'}$ be considered.  If the schemes were affine,
this would be equivalent to deciding whether to consider this tensor
product as an algebra over $\cO_X$, $\cO_Y$, $\cO_S$, etc.  Likewise,
in the general case there is no canonical choice of base scheme for
the dg scheme $W'$, and either one of $X$, $Y$, $S$, or $X\times Y$
can serve for this purpose.  For example, it is easy to see $W'$ as a
dg scheme over $X$ by resolving $\cO_Y$ by a flat commutative dg algebra
over $S$ and pulling back the resolution to $X$.  Similarly, in order
to obtain a model over $S$ resolve both $\cO_X$ and $\cO_Y$ over $S$
and tensor them over $\cO_S$.

It is essential to emphasize that in general it is not possible to
present $W'$ as a dg scheme over $W$, the underived intersection.

\paragraph
\label{subsec:jusils}
For the purpose of this article we are most interested in a model of
$W'$ whose base scheme is $X\times Y$.  To obtain such a presentation
define
\[ \cO_{W'} = \cO_{\Gamma_i} \circ \cO_{\Gamma_j} = \pi_{XY,*}
(\pi_{XS}^* \cO_{\Gamma_i} \otimes_{X\times S \times Y} \pi_{SY}^*
\cO_{\Gamma_j}), \] 
the convolution of the kernels $\cO_{\Gamma_i}\in \D(X\times S)$ and
$\cO_{\Gamma_j} \in \D(S \times Y)$.  Here $\Gamma_i \subset X \times
S$, $\Gamma_j\subset S\times Y$ are the graphs of the inclusions
$i:X\hookrightarrow S$, $j:Y\hookrightarrow S$, and $\pi_{XS}$,
$\pi_{SY}$ and $\pi_{XY}$ are the projections from $X\times S\times Y$
to $X\times S$, $S\times Y$, and $X\times Y$, respectively.  (We
omit the $\R$'s and $\Ld$'s in front of derived functors for
simplicity.)  The reader can easily supply the required equality of
tensor products of rings which shows that this definition of $W'$ is
quasi-isomorphic to the previous ones.

Note that the kernels $\cO_{\Gamma_i}$ and $\cO_{\Gamma_j}$ induce the
functors $i_*:\D(X) \ra \D(S)$ and $j^*:\D(S) \ra \D(Y)$.  Since
$\cO_{W'}$ is the convolution of these kernels, we conclude that
$\cO_{W'}$ is the kernel of the dg functor $j^* i_* : \D(X) \ra \D(Y)$.

This fact allows us to connect with our earlier discussion
in~(\ref{subsec:extss}).  Indeed, in order to guarantee the
degeneration of the spectral sequence computing $\Ext^*_S(i_*F, j_*G)$
we need to understand the functor $j^*i_*$.  Since this functor is
controlled by $W'$ as presented over $X\times Y$, this explains why we
want to understand formality properties of $W'/X\times Y$ and not over
other bases.

\paragraph
\label{subsec:replace}
There is another description of $\cO_{W'}$ as an object in $\D(X\times
Y)$ which is useful in the proof of Theorem~\ref{thm:mainthm}.  The
original problem of studying the intersection of $X$ and $Y$ into $S$
can be reformulated to study the intersection of $X\times Y$ with the
diagonal in $S\times S$.  Let $\bi$ and $\bj$ be the
embeddings of $S$ and $X\times Y$ into $S\times S$.

The derived and underived intersections in the new problem are the
same as in the old one.  The excess bundle is also the same.  However,
by replacing the original problem with the new one we have simplified
the initial situation in two ways.  First, the embedding $\bi :S
\hookrightarrow S\times S$ is now split.  Second, since the object
$\bj^* \bi_* \cO_S$ realizes $\cO_{W'}$ as an object of
$\D(X\times Y)$, the problem of understanding the functor $j^*i_*$ is
replaced by the problem of understanding the single object $\bj^*
\bi_* \cO_S$.  We have replaced the functor $j^*i_*$ by the more
complicated functor $\bj^* \bi_*$, but we only apply it to
a single object $\cO_S$ which is well behaved.

\paragraph
We now turn to questions of formality.  Given a dg scheme $Z'$ over a
scheme $Z$, with structure morphism $f:Z'\ra Z$ being affine, we shall
say that $Z'$ is {\em formal} over $Z$ if $f_* \cO_{Z'}$ is a formal
$\cO_Z$-dg algebra, that is, if there exists an isomorphism
\[ f_* \cO_{Z'} \iso \bigoplus_j \ccH^j(f_* \cO_{Z'})[-j] \]
of $\cO_Z$-dg algebras.  This is equivalent to the dg schemes $Z'$ and
$\hat{Z}'$ being isomorphic in the derived category of dg schemes over
$Z$.  (Here $\hat{Z}'$ is the dg scheme whose structure complex is
$\bigoplus_j \ccH^j(f_* \cO_{Z'})[-j]$.)  

Note in particular that this implies that the complex $f_*\cO_{Z'}$
splits in the derived category $\D(Z)$ (it is isomorphic to the sum of
its cohomology sheaves).

\paragraph
The notion of formality of a dg scheme depends on the scheme over
which we are working.  Indeed, consider a smooth subvariety $X$ of a
smooth space $S$, and let $X' = X\times_S^R X$ be the derived
self-intersection of $X$ inside $S$.  Then $X'$ is a dg scheme over
$X$ in two distinct ways (using the two projections), and hence it is
also a dg scheme over $X\times X$.  In~\cite{AriCal} the first two
authors introduced two classes,
\begin{align*}
\alpha_\HKR & \in  H^2(X, N \otimes N^\chk \otimes N^\chk)
\intertext{and}
\eta &\in H^1(X, T_X \otimes N^\chk).
\end{align*}
The results of [loc.\ cit.] and the present paper show that
\begin{itemize}
\item[--] $X'$ is formal over $X$ if and only if the HKR class
  $\alpha_\HKR$ vanishes;
\item[--] $X'$ is formal over $X\times X$ if and only if the class
  $\eta$, vanishes.
\end{itemize}
It is known~(\cite{AriCal}) that $\eta = 0$ implies $\alpha_\HKR = 0$,
but not vice-versa.  Thus $X'$ being formal over $X\times X$ implies it
is formal over $X$, but the converse can fail.
\section{The proof of the main theorem}
\label{sec:mainthm}

In this section we shall prove our main result,
Theorem~\ref{thm:mainthm}, which we restate below.  Throughout this
section we shall drop the index ``$R$'' in the notation of derived
fiber products, and write $X\times_S Y$ for the derived fiber product
$X\times_S^R Y$.  If $F$ is any vector bundle on a space $X$ we shall
let $\bbF[-1]$ denote the total space of the shifted vector bundle
$F[-1]$, i.e., the dg scheme whose structure complex is
$\S(F^\chk[1])$.   We assume that the characteristic of the ground
field $\bbk$ is either zero or large enough (in all the results below,
large enough means larger than the dimension of $S$).

\begin{Lemma}
\label{lem:spli}
Let $i:X \ra S$ be a closed embedding of smooth varieties with normal
bundle $N_{X/S}$.  A choice of splitting of $i$ (if one exists)
determines an isomorphism 
\[ X\times_S X \iso \bbN_{X/S}[-1] \]
in the derived category of dg schemes over $X\times X$, commuting with
the embeddings of $X$.
\end{Lemma}

\begin{Proof}
Let $\pi_1$ and $\pi_2$ denote the two projections from $X\times_S X$
to $X$.  We regard $X\times_S X$ as a space over $X\times X$ via the
map $(\pi_1, \pi_2)$.  Note that by the very construction of the fiber
product the compositions $i \circ \pi_1$ and $i\circ \pi_2$ are
homotopic in the $(\infty, 1)$-category of dg schemes (before deriving
it).

Fix a splitting $p$ of the embedding $i$.  Composing a homotopy
between $i\circ \pi_1$ and $i\circ \pi_2$ with $p$ we conclude that
the maps $\pi_1$ and $\pi_2$ are themselves homotopic.  Thus in the
derived category of dg schemes, where homotopic maps become equal, we
have $\pi_1 = \pi_2$. In other words the original map
\[ (\pi_1, \pi_2): X\times_S X \ra X\times X \]
is equal to the map 
\[ (\pi_1, \pi_1): X\times_S X\ra X\times X \]
and the latter obviously factors through the diagonal map: $(\pi_1,
\pi_1) = \Delta \circ \pi_1$.  Thus the structure map $(\pi_1, \pi_2)$
of $X\times_S X$ factors through $\Delta$. 

The splitting of the map $i$ implies that we are in a situation where
we can apply the main theorem of~\cite{AriCal}.  Thus (choosing one
side) there exists an isomorphism of spaces over $X$
\[ \phi: X\times_S X \stackrel{\sim}{\lra} \bbN_{X/S}[-1], \] 
where $X\times_S X$ is regarded as a space over $X$ via $\pi_1$.
Since we have seen that the structure maps from $X\times_S X$ and
$\bbN_{X/S}[-1]$ to $X\times X$ factor through the diagonal map, it follows
that $\phi$, which originally was an isomorphism over $X$, can be
regarded as an isomorphism over $X\times X$.  The compatibility with
the embeddings of $X$ is obvious from the construction of
$\phi$ in~\cite{AriCal}.
\qed
\end{Proof}

\paragraph
We now place ourselves in the context of~(\ref{subsec:diagintro}),
with $X$ and $Y$ smooth subschemes of $S$, and with $W'$ and $W$ being
their derived and underived intersections, respectively.  The maps
between these spaces are listed in the diagram below
\[ 
\xymatrix{W' \ar@/_/@{.>}[dr]_{\pi} \ar@/_1pc/[ddr]_{q'}
  \ar@/^1pc/[drrr]^{p'}& & & \\ & W\ar@/_/[ul]_-{\phi}
  \ar[rr]^p\ar[d]^-q & & X\ar[d]^-i\\ & Y\ar[rr]^j & & S.\\}
\]
The excess intersection bundle $E$ on $W$ is defined as
\[ E = \frac{T_S}{T_X+T_Y} \] 
where all the bundles above are assumed to have been restricted to
$W$.  As usual we shall denote by $\bbE[-1]$ the total space of the
shift of $E$.  
\medskip

\noindent
We begin by studying a special case of the main theorem, where the map
$i$ is split.

\begin{Proposition}
\label{easyformal}
Assume that the map $i$ is split, and fix a splitting of $i$.  Then
a choice of splitting of the short exact sequence
\[ \hspace{1.2in} 0 \ra N_{W/Y} \ra N_{X/S}|_W \ra E \ra 0 \hspace{1.15in} (*)\]
gives rise to an isomorphism of spaces over $X\times Y$
\[ \bbE[-1] \stackrel{\sim}{\lra}  X\times_S Y. \]

Conversely, existence of an isomorphism in $\D(Y)$
\[ j^* i_* \cO_X \iso q_*(\S(E^\chk[1]) \]
implies that the short exact sequence $(*)$ splits.
\end{Proposition}

\begin{Proof}
A splitting $E \ra N_{X/S}|_W$ of the short exact sequence $(*)$ gives
rise to a map $\bbE[-1] \ra \bbN_{X/S}[-1]|_W$.  We have fixed a
splitting of $i$; by Lemma~\ref{lem:spli} this gives rise to an
isomorphism over $X\times X$ 
\[ \bbN_{X/S}[-1] \stackrel{\sim}{\lra} X\times_S X, \]  
compatible with the inclusion of $X$.  Using this
isomorphism we obtain a morphism over $X\times Y$
\[ \bbE[-1] \ra \bbN_{X/S}[-1]|_W \iso (X\times_S X)\times_X W \iso
X\times_S W \ra X\times_S Y \]
where the last map comes from the inclusion $W \hookrightarrow Y$.
Checking that this morphism is an isomorphism is a local computation
which can be checked using Koszul resolutions. 

In the other direction assume that there exists an isomorphism in $\D(Y)$
\[ \phi: j^* i_* \cO_X \stackrel{\sim}{\lra} q_*(\S(E^\chk[1]). \]
Without loss of generality we can assume that $\phi$ commutes with the
natural maps of the two sides to $q_*\cO_W$.  To see this consider the
map $\phi^0$ induced by $\phi$ on $\ccH^0$ of the two complexes.  It
is an automorphism of the $\cO_Y$-module $q_*\cO_W$.  As $q_*\cO_W$ is
a quotient algebra of $\cO_Y$, $\phi^0$ is in fact an
automorphism of the ring $\cO_W$, given by multiplication by an
invertible element $s$ of $\cO_W$.  Since $\S(E^\chk[1])$ is an
$\cO_W$-algebra, multiplication by $s^{-1}$ gives
an automorphism $\psi$ of it.  The composition $q_*\psi\circ \phi$ is
a new isomorphism
\[ j^* i_* \cO_X \iso q_*(\S(E^\chk[1])) \] 
which commutes with the maps to $q_*\cO_W$, as desired.  We shall call
this new map $\phi$.

The map $\phi$ induces an isomorphism
\[ p^* i^* i_* \cO_X \iso q^* j^* i_* \cO_X\iso q^* q_* (\S(E^\chk[1])). \]
Applying $\ccH^{-1}$ to both sides gives an isomorphism
\[ N_{X/S}^\chk|_W \iso N_{W/Y}^\chk \oplus E^\chk \]
where the component $N_{W/Y}^\chk$ comes from $\ccH^{-1}$ of the
summand $q^*q_*\cO_W$ of $q^*q_*(\S(E^\chk[1]))$.  The fact that
$\phi$ is compatible with the map to $q_*\cO_W$ shows that the map
\[ N_{X/S}^\chk \ra N_{W/Y}^\chk \]
in the above decomposition is the same as the map obtained by applying
$\ccH^{-1}$ to the morphism 
\[ p^*i^*i_* \cO_X \iso q^*j^*i_*\cO_X \ra q^*q_*\cO_W. \]
A straightforward calculation with the Koszul complex shows that this
map is precisely the dual of the inclusion map 
\[ N_{W/Y} \ra N_{X/S}|_W \]
from the short exact sequence $(*)$.  Thus the direct sum
decomposition above is compatible with the maps in $(*)$, and hence
this short exact sequence must split.
\qed
\end{Proof}

\begin{Theorem}
The following statements are equivalent.
\begin{itemize}
\item[(1)] There exists an isomorphism of dg functors $\D(X) \ra \D(Y)$
\[ j^* i_* (\,-\,) \iso q_*(p^*(\, - \,) \otimes \S(E^\chk[1])) \]
\item[(2)]  There exists an isomorphism $W'\iso \bbE[-1]$ of dg
  schemes over $X\times Y$.
\item[(3)] $W'$ is formal as a dg scheme over $X\times Y$.
\item[(4)] The inclusion $W \ra W'$ splits over $X\times Y$.
\item[(5)] The inclusion $W \ra W'$ splits to first order over $X\times Y$.
\item[(6)] The short exact sequence
\[ 0 \ra T_X + T_Y \ra T_S \ra E \ra 0 \]
of vector bundles on $W$ splits.
\end{itemize}
\end{Theorem}

\begin{Proof} 
We shall prove the following chains of implications and equivalences
\begin{align*}
& (2) \Rightarrow (3) \Rightarrow (4) \Rightarrow (5) \Rightarrow (6)
\Rightarrow (2),
\intertext{and}
& (1) \Rightarrow (6) \Rightarrow (2)  \Rightarrow (1).
\end{align*}
The implications $(2) \Rightarrow (3) \Rightarrow (4) \Rightarrow (5)$
are obvious.  The implication $(5) \Rightarrow (6)$ is a dg version
of~\cite[20.5.12 (iv)]{EGA4.1} (which we restate as
Lemma~\ref{lem:splitfirstorder} below), once one notes that the
derived relative tangent bundles involved are
\begin{align*}
T_{W/W'} & = E[-2] \\
T_{W/X\times Y} & = (T_X + T_Y)[-1] \\
T_{W'/X\times Y} & = T_S[-1].
\end{align*}
The implication $(2) \Rightarrow (1)$ follows from the considerations
in~(\ref{subsec:jusils}).  Indeed, the kernel giving the functor
$j^*i_*$ is $\cO_W'$, while the kernel giving the functor
$q_*(p^*(\,-\,) \otimes \S(E^\chk[1]))$ is $\cO_{\bbE[-1]}$, and an
isomorphism between these objects gives rise to an isomorphism between
the corresponding functors.

In the other direction there is a subtle point.  Work of
To\"en~\cite{Toe} does imply that the equivalence of functors in $(1)$
guarantees an isomorphism of {\em $\cO_{X\times Y}$-modules} 
\[ \cO_{W'} \iso \cO_{\bbE[-1]}; \] 
while condition $(2)$ is the strongest requirement that the two be
isomorphic as {\em algebras}.  Hence the implication $(1) \Rightarrow
(2)$ is not automatic and will only follow indirectly from the rest of
the proof.

The implications that we still need to prove are $(6) \Rightarrow (2)$
and $(1) \Rightarrow (6)$.  We replace the initial intersection
problem with the problem of intersecting $X\times Y$ with the diagonal
in $S\times S$, as in~(\ref{subsec:replace}). We keep denoting the new
spaces and embeddings by $X$, $Y$, and $S$, $i$, $j$, etc. Thus the
new $S$ is the old $S\times S$, the new $X$ is the diagonal in the old
$S\times S$, and the new $Y$ is the old $X\times Y$.  Note that now
$i$ is split (it is the old diagonal map, hence it is split by either
projection).

We reformulate $(1)$, $(2)$, and $(6)$ of the theorem in the new
setting.  Statements $(1)$ and $(2)$ become the statements that there
exist isomorphisms
\[ j^* i_* \cO_X \iso q_* \S(E^\chk[1]) \] 
as objects of $\D(Y)$ and as commutative dg algebra objects in
$\D(Y)$, respectively.  The short exact sequence of $(6)$ becomes the
sequence
\[ 0 \ra N_{W/Y} \ra N_{X/S}|_W \ra E \ra 0. \]

We are now in the situation of Proposition~\ref{easyformal}: indeed,
the main property we need is that the map $i$ splits, and this is true
because now $i$ is the old diagonal map.  The conclusions of this
proposition are exactly the implications $(6) \Rightarrow (2)$ and
$(1) \Rightarrow (6)$ that we still needed to prove.
\qed
\end{Proof}

\begin{Lemma}
\label{lem:splitfirstorder}
Let $i: X \hookrightarrow Y$ be a closed embedding of dg schemes over
a fixed scheme $S$.  Then $i$ is split to first order over $S$ if and
only if the natural map 
\[ \bbT_{X/Y} \ra \bbT_{X/S} \]
is the zero map, where $\bbT$ denotes the tangent complex of the
corresponding morphism.
\end{Lemma}

\begin{Proof}
The proof is nothing but a restating in dg language of~\cite[20.5.12
(iv)]{EGA4.1}.
\end{Proof}
\section{Applications to orbifolds}

In this section we discuss how Theorem~\ref{thm:mainthm} can be used
to understand the structure of derived fixed loci.  In turn this allows
us to understand the structure of the free loop spaces of orbifolds.

\paragraph
\label{subsec:fixedloc}
We review the setup in~(\ref{subsec:derfixloc}).  Let $Z$ be a smooth
variety over a field $\bbk$, and let $\phi$ be an automorphism of $Z$ of
finite order $n$. Since $\phi$ is of finite order its fixed locus,
which we shall denote by $W$, is scheme-theoretically smooth.  We shall
assume that the characteristic of $\bbk$ is either zero or greater than
$\max(n, \dim Z)$.

Note that the ordinary fixed locus $W$ can be understood as the
intersection
\[ W= \Delta \times_{Z\times Z} \Delta^\phi, \] 
where $\Delta$ and $\Delta^\phi$ denote the diagonal in $Z\times Z$
and the graph of $\phi$, respectively.  As such the expected dimension
of $W$ is zero.  Whenever $\dim W>0$ it is important to understand the
failure of the spaces in this intersection problem to meet
transversally, by studying the derived intersection space
\[ W' = \Delta\times_{Z\times Z}^R \Delta^\phi. \] 
We shall sometimes call this space the derived fixed locus of $\phi$.

\paragraph
Theorem~\ref{thm:mainthm} shows that in order to understand the
structure of $W'$ we need to study the short exact sequence 
\[ 0 \ra T_\Delta + T_{\Delta^\phi} \ra T_{Z\times Z} \ra E \ra 0 \]
where $E$ is the excess bundle for this intersection problem.  We
shall prove that this sequence is always split, under the assumptions
we made for the characteristic of $\bbk$.  We begin with a lemma.

\begin{Lemma}
\label{two:seq}
In the setup of Theorem~\ref{thm:mainthm}, assume that the map $X\ra
S$ is split to first order. Then the short exact sequence
\[ 0 \ra N_{W/Y} \ra N_{X/S}|_W \ra E \ra 0 \]
splits if and only if the six equivalent statements of
Theorem~\ref{thm:mainthm} are all true.
\end{Lemma}

\begin{Proof}
It is easy to see that the two conditions of the lemma imply that the
short exact sequence of $(6)$ of Theorem~\ref{thm:mainthm} splits.
Equivalently, these two conditions are what was used in the proof of
Theorem~\ref{thm:mainthm} after changing the problem to an
intersection of the diagonal with $X\times Y$.
\qed
\end{Proof}

\begin{Theorem}
\label{der:fix:loc}
Assume we are in the setup of~(\ref{subsec:fixedloc}). Then the
derived fixed locus $W'$ is isomorphic, as a dg scheme over $Z\times
Z$, to the total space over $W$ of the dg vector bundle
$T_W[-1]$,
\[ W' \iso \bbT_W[-1]. \]
\end{Theorem}
\vspace{-4mm}

\begin{Proof}
It is easy to see that the excess intersection bundle $E$ for this
intersection problem is $(T_Z)_\phi$, the bundle on $W$ of
coinvariants of the action of $\phi$ on $T_Z$,
\[ (T_Z)_\phi= \frac{T_Z}{\langle v - \phi(v) \rangle}. \] 
We now apply Lemma \ref{two:seq}. The embedding $Z\ra Z\times Z$ is
split to first order (it is actually split). The map $N_{Z/Z\times
  Z}\ra E$ is given by the natural projection
\[ T_Z\ra \frac{T_Z}{\langle v - \phi(v) \rangle}. \] 
If the characteristic of $\bbk$ is 0 or prime to $n$, the averaging
map $(T_Z)_\phi\ra T_Z$ given by  
\[ t\mapsto \frac{1}{n}\sum_{i=1}^n \phi^i(t) \]
splits the projection above.  Finally, with the same assumptions on
the characteristic of $\bbk$, the bundles of invariants and
coinvariants are naturally isomorphic: $(T_Z)_\phi \iso (T_Z)^\phi$
and the latter is precisely $T_W$.
\qed
\end{Proof}

\paragraph
\label{sec:orbifolds:setup}
We apply the above theorem to the study of orbifolds. Let $G$ be a
finite group acting on a smooth variety $Z$, and denote the
quotient stack $[Z/G]$ by $\cZ$. We are interested in understanding
the relationship between the inertia stack of $\cZ$,
\[I\cZ=\cZ\times_{\cZ\times \cZ}\cZ,\]
and the free loop space of the corresponding derived intersection,
\[L\cZ=\cZ\times^R_{\cZ\times \cZ}\cZ.\]
We organize these spaces and the maps between them in the diagram below:
\[ 
\xymatrix{L\cZ \ar@/_/@{.>}[dr]_{\pi} \ar@/_1pc/[ddr]_{q'}
  \ar@/^1pc/[drrr]^{p'}& & & \\ & I\cZ\ar@/_/[ul]_-{\phi}
  \ar[rr]^p\ar[d]^-q & & \cZ\ar[d]^-\Delta\\ & \cZ\ar[rr]^\Delta & & \cZ\times \cZ.\\}
\]

\paragraph
We wish to realize the above diagram of (dg) stacks as the global
quotient by the fixed group $G\times G$ of a similar diagram of (dg)
schemes.  This allows us to reduce the problem of understanding the dg
stack $L\cZ$ and its maps to $I\cZ$ and $\cZ$ to the parallel
problem of understanding the corresponding dg schemes and maps. 

As originally formulated the diagonal map $\Delta$ is a map between
the global quotient stacks $\cZ = [Z/G]$ and $\cZ\times \cZ =
[(Z\times Z)/(G\times G)]$.  We wish to replace the presentation
$[Z/G]$ of $\cZ$ by a different presentation of the same stack, but
where the group we quotient by is $G\times G$.  Consider the action of
$G\times G$ on $Z\times G$ given by
\[(h,k).(z,g)\mapsto (h.z,kgh^{-1}).\] 
Note that the second copy of $G$ acts freely on $Z\times G$, thus yielding
an isomorphism
\[[(Z\times G)/(G\times G)]\iso [Z/G].\]
With this presentation the diagonal map $\cZ\ra \cZ\times
\cZ$ becomes the quotient by $G\times G$ of the equivariant map of
spaces 
\begin{align*}
\bar{\Delta} :Z\times G& \ra Z\times Z \\
(z,g) & \mapsto (z,g.z). 
\end{align*}
Summarizing the above discussion, the main diagram
in~(\ref{sec:orbifolds:setup}) is obtained by taking the quotient by
$G\times G$ of the spaces and maps in the diagram below:
\[ 
\xymatrix{LZ \ar@/_/@{.>}[dr]_{\bar{\pi}} \ar@/_1pc/[ddr]_{\bar{q}'}
  \ar@/^1pc/[drrr]^{\bar{p}'}& & & \\ & IZ\ar@/_/[ul]_-{\bar{\phi}}
  \ar[rr]^{\bar{p}}\ar[d]^-{\bar{q}} & & Z\times G\ar[d]^-{\bar{\Delta}}\\ & Z\times
  G\ar[rr]^{\bar{\Delta}} & & Z\times Z.\\}
\]
Here we have denoted by $LZ$ and $IZ$ the corresponding derived and
underived fiber products, respectively.

\paragraph
\label{subsec:quotG}
For $g\in G$ denote by $\Delta^g$ the subvariety of $Z\times Z$ which
is the graph of the action of $g$ on $Z$,
\[ \Delta^g = \{ (z, g.z)~|~z\in Z\}. \] 
The space $IZ$ decomposes as the disjoint union
\[ IZ = \coprod_{g, h\in G} \Delta^g \cap \Delta^h, \] 
and similarly for $LZ$ (where the intersection is replaced by the
derived intersection).  Since the second copy of $G$ acts freely on
$X\times G$ it follows that it also acts freely on $IZ$ and $LZ$
(which can be thought of as subvarieties of $X\times G$).  This allows
us to further simplify the calculation of $I\cZ$ and $L\cZ$ by first
taking the quotient of $IZ$ and $LZ$ by the second copy of $G$ (which
still yields a space), leaving the first copy of $G$ to quotient by
later.  This amounts to replacing in the above calculation of a
(derived) intersection the horizontal map $\bar{\Delta}$ by just the
map $\Delta:Z \ra Z\times Z$, while the vertical map $\bar{\Delta}$
stays the same.  We shall abuse notation and denote the new derived and underived
intersection spaces by the old names of $LZ$ and $IZ$.  They fit in
the diagram
\[ 
\xymatrix{LZ \ar@/_/@{.>}[dr]_{\bar{\pi}} \ar@/_1pc/[ddr]_{\bar{q}'}
  \ar@/^1pc/[drrr]^{\bar{p}'}& & & \\ & IZ\ar@/_/[ul]_-{\bar{\phi}}
  \ar[rr]^{\bar{p}}\ar[d]^-{\bar{q}} & & Z\times G\ar[d]^-{\bar{\Delta}}\\ & Z\ar[rr]^{\Delta} & & Z\times Z.\\}
\]

\paragraph
Observe that after this reduction the space $IZ$ is just the disjoint
union 
\[ IZ = \coprod_{g\in G} Z^g, \]
while $LZ$ has the same decomposition as a disjoint union, but the
fixed loci $Z^g$ are replaced by their derived analogues $(Z^g)'$.  The action
of $h\in G$ shuffles these fixed loci by sending $Z^g$ to
$Z^{hgh^{-1}}$.
\medskip

\noindent
Applying Theorem~\ref{der:fix:loc} immediately yields
Theorem~\ref{thm:orb1}, which we restate below.

\begin{Theorem}
\label{thm:orb1new}
Let $\cZ$ be a global quotient orbifold.  Then there exists a
canonical isomorphism 
\[ \exp: \bbT_{I\cZ}[-1] \stackrel{\sim}{\lra} L\cZ \]
between the shifted tangent bundle $\bbT_{I\cZ}[-1]$ of its inertia
orbifold and its free loop space $L\cZ$.
\end{Theorem}

\begin{Proof}
Theorem~\ref{der:fix:loc} implies that the derived fixed loci $(Z^g)'$
are isomorphic to the total spaces $(\bbT_Z)_g[-1]$, and it is
immediate to see that $(T_Z)_g \iso T_{Z^g}$.  Thus we
have an isomorphism
\[ \bbT_{IZ}[-1] \stackrel{\sim}{\lra} LZ; \]
the isomorphism in the statement of the theorem is nothing but the
quotient of this isomorphism by the action of $G$.
\qed
\end{Proof}

\paragraph
Corollary~\ref{cor:orb2} follows easily from the above theorem once
one remembers that 
\begin{align*}
HH_*(\cZ) & = \R \Gamma(\cZ, \Delta^* \Delta_* \cO_{\cZ}) =
\R\Gamma(Z, \bar{q}'_*\cO_{LZ})^G \\
\intertext{and}
HH^*(\cZ) & = \R\Hom^*_{\cZ\times \cZ}(\Delta_* \cO_\cZ,
\Delta_*\cO_\cZ) = \R\Gamma(\cZ, \Delta^!\Delta_* \cO_\cZ) \\
& = \R\Gamma(Z, (\bar{q}'_* \cO_{LZ})^\chk)^G.
\end{align*}
Here $\bar{q}'$ is the map $LZ \ra Z$ from~(\ref{subsec:quotG}).

\paragraph
Theorem~\ref{thm:orb1new} also highlights a somewhat striking
difference between the behavior of disconnected Lie groups in derived
algebraic geometry versus classical geometry.

We think of the free loop space $L\cZ$ of an orbifold $\cZ$ as a
family of (homotopy) groups $\Omega_z \cZ$ parametrized by $z\in\cZ$.
These groups are not in general connected, having components indexed
by the inertia group of $z$.  Our theorem shows that the exponential
map we have defined is an isomorphism between an appropriate number of
copies of the Lie algebra of the group $\Omega_z\cZ$ and the group
itself.

This situation is different when one studies disconnected Lie groups
in the classical setting.  Let $H$ be a disconnected Lie group, let
$H^0$ be the connected component of the identity (a normal subgroup of
$H$), let $\gh$ denote the Lie algebra of $H$ (i.e., the Lie algebra
of $H^0$), and let $G = H/H^0$ denote the group of components.
Suppose the homomorphism $H\ra H/H_0=G$ admits a right inverse. Once
the inverse is fixed, we can view $G$ as a subgroup of $H$, and $H$
decomposes as a semi-direct product of $H_0$ and $G$.

In this setting there is no natural exponential map that covers {\em
  all} of $H$ (or at least a formal neighborhood of $G$ in $H$).
Indeed, one could translate the usual exponential map around the
origin of $H$ to the other components of $H$, but this involves a
choice -- whether to use left or right translation by elements of $G$.
The only natural map to a neighborhood of $g\in G$ is the restriction
of the exponential map to $(H^0)^g$, which will map into the
centralizer of $g$.  This map will be far from surjective, unlike the
derived case where it is an isomorphism.

\end{document}